\documentclass[12pt]{article}
\usepackage{amsmath}
\usepackage{amssymb}
\def \tp {{\rm tp}}

\def \dcl {{\rm dcl}}

\def \non-ort {{\not\!\!\bot}}

\def \proof {\noindent{\it Proof.\ }}
\def \qed{\hfill\rule{4pt}{9pt} \newline\medbreak}

\def\proclaim #1. #2\par{\medbreak
  \noindent{\bf#1 \enspace}{\sl#2}\par
  \ifdim\lastskip<\medskipamount \removelastskip\penalty55\medskip\fi}

\newtheorem{theorem}{Theorem}
\newtheorem{definition}[theorem]{Definition}
\newtheorem{proposition}[theorem]{Proposition}

\newtheorem{corollary}[theorem]{Corollary}

\newtheorem{remark}[theorem]{Remark}

\newtheorem{claim}{Claim}[theorem]

\title{Definability of initial segments}
\author{
Saharon S{\sc helah}\footnote{
The first author was partially supported by the United
States-Israel Binational Science Foundation.
This work was done while the first author was visiting 
Japan in summer 2000.
He expresses his thanks to the Japanese Association
of Mathematical Sciences for inviting him. 
The number of this article in Shelah's publication list  is 767.
}  \
 and Akito T{\sc suboi}\footnote{
The second author was partially supported by
Japanese Grant-in-Aid for Scientific Research (C2) 11640100.
}
}
\date{}
\begin{document}
\maketitle

\section{Introduction}

Let $T$ be a first order theory formulated in the language $L$ and
$P$ a new relation symbol not in  $L$.
Let $\varphi(P)$ be an $L \cup \{P\}$-sentence.
Let us say that $\varphi(P)$ defines $P$ implicitly in $T$
if $T$ proves $\varphi(P) \wedge \varphi(P') \rightarrow
\forall x(P(x) \leftrightarrow P(x'))$.
Beth's definability theorem states that if
$\varphi(P)$ defines $P$ implicitlly in $T$ then
$P(x)$ is equivalent to an $L$-formula.

However, if we consider implicit definability in a
given model alone, the situation changes.
For a more precise explanation,
let us say that a subset $A$ of a given model $M$
of $T$ is implicitly definable  if there exists a sentence $\varphi(P)$
such that $A$ is the unique set with $(M,A) \models \varphi(P)$.
It is easy to find a structure
in which two kinds of definability (implicit definability and first order definability) are different.
For example, let us consider the structure $M=(\mathbb N \cup \mathbb Z,<)$,
where $<$ is a total order such that any element in the
$\mathbb Z$-part is greater than any element in the
$\mathbb N$-part.
The $\mathbb N$-part is not first order definable in $M$, because
the theory of $M$ admits quantifier elimination after adding the
constant $0$ (the least element) and
the successor function to the language.
But the $\mathbb N$-part  is implicitly definable in $M$,
because it is the unique non trivial initial segment
without the last element.
On the other hand, for a given structure,
we can  easily  find an elementary extension in which
two notions coincide.

In this paper, we shall consider implicit definability of the
standard part $\{0,1,...\}$  in  nonstandard models of Peano
arithmetic ($PA$).
It is needless to say that the standard part of a nonstandard
model of $PA$
is not first order definable.  As is stated above, there is a model
in which every set defined implicitly is first order definable.
So we ask whether there is a model of $PA$ in which the
standard part is implicitly definable.

In  \S 1,  we define a certain class of formulas, and show that
in any model of $PA$ the standard part
is not implicitly defined by using such formulas.

\S 2 is the main section of the present paper,
we shall construct a model of $PA$ in which the standard
part is implicitly defined.
To construct such a model, first we assume a set theoretic
hypothesis $\diamondsuit_{S_\lambda^{\lambda^+}}$,
which is an assertion of  the existence of a very general set.
Then we shall eliminate the hypothesis using absoluteness
for the existence of a model having a tree structure with a certain
property.

In this paper $L$ is a first order countable language.
$L$-structures are denoted by $M$, $N$, $M_i$, $\cdots$.
We do not strictly distinguish  a structure and its universe.
$A$, $B$, $\cdots$ will be used for denoting subsets of
of some $L$-structure.  Finite tuples of elements from
some $L$-structure are denoted by $\bar a$, $\bar b, \cdots$.
We simply write $A \subset M$ for expressing that $A$ is
a subset of the universe of $M$.

\section{Undefinability result}

Let us first recall the definition of implicit definability.

\begin{definition}\rm
Let $M$ be an $L$-structure.  Let $P$ be a unary
second order variable. A subset $A$ of $M$ is said to be
{\it implicitly definable} in $M$ if there is an
$L \cup \{P\}$-sentence $\varphi(P)$ with parameters
such that $A$ is the unique solution to $\varphi(P)$, i.e.
$\{A\}=\{B \subset M: M \models \varphi(B)\}$.
\end{definition}

In this section $L$ is the language $\{0,1,+,\cdot,<\}$, and
$PA$ denotes the Peano arithmetic formulated in $L$.
We shall prove that the standard part  is not
implicitly definable in any model of $PA$
by using a certain form of formulas.
We fix a model $M$ of $PA$, and work on $M$.

\begin{definition}\rm
An $L \cup \{P\}$-formula $\varphi(\bar y)$ (with parameters) will
be called simple if it is equivalent to a prenex normal form
$
Q_1 \bar x_1\cdots Q_n \bar x_n
[
P(f(\bar x_1,...,\bar x_n, \bar y)) \rightarrow
P(g(\bar x_1,...,\bar x_n,\bar y))
]
$
where
$Q_i$'s are quantifiers and
$f$ and $g$ are definable functions.
If $Q_1=\forall$  then $\varphi$ will be called a simple
$\Pi_n$-formula. Similarly it is called a simple $\Sigma_n$-formula
if $Q_1=\exists$.
\end{definition}

\begin{remark}\rm
If $P$ is an initial segment of $M$, then
\begin{enumerate}
\item $a_1 \in P \wedge a_2 \in P$ is equivalent to
$\max\{a_1,a_2\} \in P$;
\item $a_1 \in P \vee a_2 \in P$ is equivalent to
$\min\{a_1,a_2\}\in P$.
\end{enumerate}
An $L$-formula $\varphi(\bar x)$ is equivalent to a formula of
the form $P(f(\bar x))$, where $f$ is a definable function
such that $f(\bar x)=0$ if $\varphi(\bar x)$ holds and
$f(\bar x)=a$ ($a$ is a nonstandard element) otherwise.
In what follows,
an initial segment $I \subsetneq M$ will be called a {\it cut }
if $I$ is closed under successor. The statement that $P$ is a
cut is expressed by a simple $\Pi_2$-formula.
\end{remark}

We shall prove that
the standard part is not implicitly definable
by a finite number of simple $\Pi_2$-formulas.
In fact we can prove more.

\begin{proposition}\label{proposition}
Let $I_0$ be a cut of $M$ with $I_0 <a$ i.e.
any element of $I_0$ is smaller than $a$.
Let $\{\varphi_i(P):i\leq n\}$ be a finite set of simple formulas.
If  $I_0$ satisfies $\{\varphi_i(P):i\leq n\}$, then
there is another cut  $I<a$ which also satisfies
$\{\varphi_i(P):i\leq n\}$.
\end{proposition}

Let us say that a cut $I$ is approximated by a decreasing
$\omega$-sequence, if
there is a definable function $f(x)$ with
$
I=\{a \in M: \ (\forall m \in \omega)\ a \leq f(m) \}.
$
Similarly we say that $I$ is approximated by an increasing
$\omega$-sequence if there is a definable  function
$g(x)$ with $
I=\{a \in M: (\exists m \in \omega)\ a \leq g(m)  \}.
$
Notice that no cut of $M$ is approximated by both a decreasing
$\omega$-sequence
and an increasing $\omega$-sequence.

\bigbreak\noindent
{\it Proof of Proposition \ref{proposition}:}
For $i \leq n$, let $\varphi_i(P)$ have the form
$
\forall \bar x \exists \bar y
 [
P(f_i(\bar x,\bar y))  \rightarrow P(g_i(\bar x,\bar y))
].
$
By the remark just after Proposition \ref{proposition}, we can assume
that $I_0$ cannot be approximated by a decreasing
$\omega$-sequence.
We shall show that there is an initial segment $I$ with
$I_0 \subsetneq I < a$ and
$M \models \bigwedge_{i \leq n}\varphi_i(I)$.
Since $I_0$ satisfies $\varphi_i(P)$, for each
$b_0 \in M$ with $I_0<b_0<a$, we have
$
M \models
\bigwedge_{i \leq n}
\forall \bar x \exists \bar y
[
f_i(\bar x, \bar y) \in \omega \rightarrow
g_i(\bar x, \bar y) \leq b_0
]
$.
By overspill there is an element $b_1$ with $I_0 < b_1 <b_0$
such that
$$
M \models
\bigwedge_{i \leq n}
\forall \bar x \exists \bar y[
f_i(\bar x, \bar y) \leq b_1 \rightarrow
g_i(\bar x, \bar y) \leq b_0
]
$$
By choosing maximum such $b_1 < b_0$, we may assume
that $b_1 \in \dcl(\bar a,b_0)$, where
$\bar a$ are parameters necessary for defining $f_i$'s and $g_i$'s.
So we can choose  an $L(\bar a)$-definable function, $h(x)$
such that
(i) $I_0 <b <a$ implies $I_0<h(b)<b$ and (ii)
$M \models
\bigwedge_{i \leq n}
\forall x \exists y
[
f_i(\bar x, \bar y) \leq h(b) \rightarrow
g_i(\bar x, \bar y) \leq b
]$, for any nonstandard $b \in M$.

By using recursion we can choose a definable function
$l(x)$ with
$
l(m)=h^m (a)
$ (the $m$-time application of $h$)
 for each $m \in \omega$.
Now we put
$$
I=
\{
d \in M:  ( \forall m \in \omega)\ d \leq l(m)
\}.
$$
Since $m < h(m)$ holds for any $m \in \omega$,
by overspill, there is a nostandard $m^*$ such that
$m^* < h(m^*)$. This shows that
 $I$ is an initial segment different from $I_0$.
 Now we show:
\proclaim Claim. For all $i \leq n$ and for all $\bar d \in M$,
there is $\bar e \in M$ such that
$$
f_i(\bar d, \bar e) \in I\ \rightarrow \ g_i(\bar d, \bar e) \in I.
$$
\par

Let  $d \in M$ and $i \leq n$ be given.  We can assume that
$
\forall y  (f_i(\bar d, \bar y) \in I)
$
holds  in $M$.
So by the definitions of $I$ and $l$, for all $k \in \omega$,
we have
$
M \models
\forall y   (f_i(\bar d, \bar y) \leq l(k)).
$
Hence, for some nonstandard  $k^* \in M$ with
$k^* \leq l(k^*)$,  we have
$$
M \models
\forall \bar y (f_i(\bar d, \bar y)
\leq l(k^*)).
$$
On the other hand, by our choice of $h$ and $l$, we can
find $\bar e$ with
$$
M \models
f_i(\bar d, \bar e) \leq l(k^*) \rightarrow
g_i(\bar d, \bar e) \leq l(k^*-1).
$$
Hence, for this $\bar e$, we have
$
g_i(\bar d, \bar e) \leq l(k^*-1) \in I.
$ \qed

\begin{corollary}
The standard part is not implicitly definable by a finite number of
simple $\Sigma_3$-formulas.
\end{corollary}

\section{Definability result}

In this section we aim to prove the following theorem:
\begin{theorem}\label{theorem}
There is a model of $PA$ in which the standard part is
implicitly definable.
\end{theorem}
Instead of proving the theorem, we prove a more general result
(Theorem \ref{main}),  from which Theorem \ref{theorem} easily follows.
For stating the result, we need some preparations.

We assume the language $L$ contains a binary predicate
symbol $<$, a constant symbol $0$ and a unary function symbol $S$.
We fix a complete $L$-theory $T$ with
a partial definable function $F(x,y)$ such that the following sentences
are members of  $T$:
\begin{itemize}
\item $<$ is a linear order with the first element $0$;
\item For each $x$, $S(x)$ is the immediate successor of $x$ with respect to $<$;
\item $\forall y_1,...,y_n \forall z_1,...,z_n\exists x
(\bigwedge_{i\neq j} y_i \neq y_j \ \rightarrow \
\bigwedge_{i=1}^{n}F(x,y_i)=z_i)$\ \  (for $n \in \omega$).
\end{itemize}

\begin{remark}\rm
Any completion of $PA$ satisfies our requirements stated
above.
\end{remark}
Let $P$ be a new unary predicate symbol not in $L$.
Throughout this section $\psi^*(P)$ is
the conjunction of the following $L\cup \{P\}$-sentences:
\begin{enumerate}
\item $P$ is a cut
(non-empty proper initial segment closed under $S$), i.e.

$\neg (\forall x P(x)) \wedge P(0) \wedge
\forall x \forall y (P(y) \wedge x<y
\rightarrow P(x)) \wedge \forall x (P(x) \rightarrow P(S(x)))$;
\item For no $x$ and $z$ with $P(z)$, is $\{F(x,y):y<z \} \cap P$
 unbounded in $P$, i.e.
$\forall x \forall z [ P(z) \rightarrow \exists w
( P(w) \wedge \forall y (P(F(x,y)) \rightarrow F(x,y) < w))]$.
\end{enumerate}
It is clear that in any model $M$ of $T$, the ``standard'' part
$I=\{S^n(0):n \in \omega\}$ satisfies $\psi^*(P)$, i.e.
the sentence $\psi^*(P)$ holds in the $L \cup \{P\}$-structure $(M,I)$.

\begin{definition}\rm
A model $M$ of $T$ will be called $\psi^*$-appropriate if
the following two conditions are satisfied:
\begin{enumerate}
\item $M \neq \{S^n(0):n \in \omega\}$;
\item If $(M, I) \models \psi^*(P)$ then
(a) $I=\{S^n(0):n \in \omega\}$ or (b) $I$ is definable
in $M$ by an $L$-formula with parameters.
\end{enumerate}
\end{definition}

\begin{remark}\rm
In case that $T$ is a completion of $PA$, the part (b)
of the condition 2 in the above definition does not occur, because
in any model of $T$ no definable proper subset  is closed under $S$.
\end{remark}

\begin{theorem}\label{main}
There is an appropriate model of $T^*$.
\end{theorem}

We shall prove the theorem above by a series of claims.
For a period of time, we fix an infinite cardinal $\lambda$.
First we need some definition.

\begin{definition}\rm
Let $M$ be a model of $T$ and $\varphi(x, \bar a)$ a formula
with  parameters  from $M$. We say that $\varphi(x,  \bar a)$ is
$\Gamma^{\rm sind}_F$-big (in $M$) if in some (any)
$|T|^+$-saturated model $N \succ M$ there is $A \subset  N $
with $|A| \leq |T|$
such that for any finite number of distinct elements
$a_1,...,a_n \in  N  \setminus A$,
and any elements $b_1,...,b_n \in  N $, we have
$$
N \models
\exists x [\varphi(x,\bar a) \wedge \bigwedge_{i=1}^nF(x,a_i) =b_i].
$$
In the above definition, if $\lambda=\aleph_0$, we replace the
condition $|A| \leq |T|$ by $|A|<\aleph_0$.
\end{definition}
Let us briefly recall the definition of {\it bigness} defined in [2].
Let $R \notin L$ be a unary predicate symbol.
A statement (or an infinitary $L \cup \{R\}$-sentence) $\Gamma(R)$
 is called a notion of
bigness for $T$,  if any model $M$ of $T$ satisfies the following axioms,
for all formulas $\varphi(x,\bar y)$ and $\psi(x,\bar y)$
(where $\Gamma(\varphi(x,\bar y))$ means that setting
$R(x)=\varphi(x,\bar y)$ [so $\bar y$ is a parameter] makes
$\Gamma$ true):
\begin{enumerate}
\item $\forall \bar y (\forall x(\varphi \rightarrow \psi) \wedge
\Gamma(\varphi) \ \rightarrow \ \Gamma(\psi))$;
\item $\forall \bar y (\Gamma(\varphi \vee \psi) \ \rightarrow \
\Gamma(\varphi) \vee \Gamma(\psi))$;
\item $\forall \bar y (\Gamma(\varphi) \rightarrow \exists^{ \geq 2} x
\varphi)$;
\item $\Gamma(x=x)$.
\end{enumerate}
Now let $\Gamma(\varphi)$ be the statement
``$\varphi$ is $\Gamma_F^{\rm sind}$-big''.  Then this $\Gamma$  satisfies the above
four axioms:  It is easy to see that our $\Gamma$ saitsfies
Axioms 1, 3 and 4.  So let us prove Axiom 2.
Suppose that neither $\varphi$ nor $\psi$ are  big.
Let $M$ be a model of $T$ and $N \succ M$ be $|T|^+$-satrurated.
Let $A$ be a subset of $ N $ of cardinality $\leq |T|$.
Since $\varphi$ is not big, $A$ cannot witness the definition of
bigness,  so there are
a finite number of elements
$a_1,...,a_n \in  N  \setminus A$ with no repetition
and $b_1,...,b_n \in  N$
such that
$N \models \forall x [\bigwedge_{i \leq n} F(x,a_i)=b_i \rightarrow
\neg \varphi(x)]$.
Since $\psi$ is not big,
$A'=A \cup \{a_1,...,a_n\}$ cannot witness the definition of
bigness, hence there are $a_{n+1},...,a_m \in
 N  \setminus (A \cup \{a_1,...,a_n\})$ with no repetition and
$b_{n+1},...,b_m \in N$ such that
$N \models \forall x [\bigwedge_{n+1 \leq i \leq m}
F(x,a_i)=b_i \rightarrow \neg \psi(x)]$.  So $N \models
\forall x [\bigwedge_{i \leq m} F(x,a_i)=b_i \rightarrow
\neg (\varphi(x) \vee \psi(x))]$.  Since $A$ was chosen
arbitrarily, this shows  that $\varphi \vee \psi$ is not big.

\medbreak
For simplicity we assume $\lambda > |T|$. (This assumption
is for simplicity only.)

\begin{claim}
(Under $\diamondsuit_{S_\lambda^{\lambda^+}}+\diamondsuit_\lambda$,
where $\lambda=\lambda^{<\lambda}$, $S_\lambda^{\lambda^+}=\{\delta<\lambda^+:{\rm cf}(\delta)=\lambda\}$)
There are a continuous elementary chain
$\langle M_i:i<\lambda^+ \rangle$ of models of  $T$ and
a sequence $\langle a_i: i< \lambda^+ \rangle$ of elements
$a_i \in M_{i+1} \setminus M_i$  such that
\begin{description}
\item{(a)} $|M_i|= \lambda$;
\item{(b)} $M_i$ is saturated except when $\aleph_0 \leq {\rm cf}(i)
\leq \lambda$;
\item{(c)} $\tp_{M_{i+1}}(a_i/M_i)$
is $\Gamma_F^{\rm sind}$-big, i.e. each formula in it is $\Gamma_F^{\rm sind}$-big.
\item{(d)} $M_i \subset \{F^{M_{i+1}}(a_i,c): M_i \models c<b\}$
if $b \in M_i \setminus \{S^n(0): n \in \omega\}$,
\item {(e)} if $(C_1,C_2)$ is a Dedekind cut of
$M=\bigcup_{i < \lambda^+}M_i$ of cofinality
$(\lambda^+,\lambda^+)$ then $C_1$ is
a subset of $M$ definable with parameters.
(A Dedekind cut of $M$ of cofinality $(\mu_1,\mu_2)$ is a pair
$(C_1,C_2)$ such that
(i) $M=C_1 \cup C_2$,
(ii) $\forall x \in C_1 \forall y \in C_2 [x<^M y]$,
(iii) the cofinality of $C_1$  with respect to $<$ is
$\mu_1$ and
(iv) the coinitiality of $C_2$ (i.e. the cofinality of $C_2$
with respect to the reverse ordering)  is $\mu_2$.
\end{description}
\end{claim}

\proof  See [2]. For more details, see [3].

\bigbreak
Now we expand the language $L$ by adding new binary predicate
symbols.
Let $L^* = L \cup \{ E_1,E_2,<_{\mathrm les}, <_{\mathrm tr} \}$.
We expand the $L$-structure $M$ defined in claim A
to an $L^*$-structure $M^*$ by the following interpretation.
For $a \in M$, let $i(a)= \min \{ i < \lambda^+: a \in M_{i+1}\}$.
\begin{enumerate}
\item $E_1^{M^*}=\{(a,b): i(a)=i(b)\}$;
\item $E_2^{M^*}=\{(a,b): i(a)=i(b) \text{ and
$M \models (c<a \equiv c<b)$ for every
$c \in M_{i(a)}$} \}$, \\
In other words, $(a,b) \in E_2^{M^*}$ iff $a$ and $b$ realize
the same Dedekind cut of $M_{i(a)} (=M_{i(b)})$;
\item $<_{\mathrm lev }^{M^*}=\{(a,b): i(a)<i(b)\}$;
\item $<_{\mathrm tr }^{M^*}=\{(a,b): i(a)<i(b) \text{ and
$M \models (c<a \equiv c<b)$ for every
$c \in M_{i(a)}$} \}$.
\end{enumerate}

The relation $<_{\mathrm tr}$ defines a preorder on $M^*$ and
induces a tree structure on the $E_2$-equivalence classes.
This tree structure $(M^*/E_2, <_{\mathrm tr})$ is a definable object of
${M^*}^{\rm eq}$.
(We do not use a new symbol for the order induced by
$<_{\mathrm tr}$.)
Simiarly $<_{\mathrm lev}$  induces a linear order on
the $E_1$-equivalence classes.
Let  $R$ be the definable function which maps $a_{E_2}$ to
$a_{E_1}$. $R$ is considered as a rank function which
assigns a level to each node of the tree.
Then $\langle <_{\mathrm tr}, <_{\mathrm lev},R \rangle$ is
an {\it $L^*$-tree}  in the sense of [1].
A subset $B$ of $M^*/E_2$ will be called a {\it branch} of the tree if
(i) it is linearly ordered by $<_{\mathrm tr}$, (ii) $a_{E_2} \in B$ and
$b \leq_{\mathrm tr} a$ imply $b_{E_2} \in B$ and (iii)
the set $\{R(a_{E_2}):a_{E_2} \in B\}$ of all levels in $B$ is
unbounded in $M^*/E_1$.

\begin{claim}
Every branch of
the tree $(M^*/E_2,<_{\mathrm tr}, <_{\mathrm lev},R)$  is definable in $M^*$.
\end{claim}

\proof
Let  $B$ be a branch of the tree
$(M^*/E_2,<_{\mathrm tr}, <_{\mathrm lev},R)$.
We show that $B$ is definable in $M^*$.
Let $I$ be the $<$-initial segment determined by $B$,
i.e.
$$
I=
\{a \in M^*: M^* \models
(\forall b_{E_2} \in B)(\exists c_{E_2}\in B)[b_{E_2} <_{\mathrm tr}
c_{E_2} \wedge a < c ]\}.
$$
It is easy to see that $I$ and $B$ are interdefinable in $M^*$.
In fact, we have $b_{E_2} \in B$ if and only if
there exist $c \in I$ and $d \in M^* \setminus I$ such that
\begin{itemize}
\item $b_{E_2}$ intersects the
interval $[c,d]$,
\item  if $b_{E_2} \subset I$ then
any other $ b'_{E_2} $ with $[c,d] \cap I \cap b'_{E_2} \neq \emptyset$
has a strictly larger level than $b_{E_2}$ and
\item
if $b_{E_2} \subset M^* \setminus  I$ then any other
$ b'_{E_2} $ with $[c,d] \cap (M^* \setminus I) \cap b'_{E_2}
\neq \emptyset$ has
a strictly larger level than $b_{E_2}$.
\end{itemize}

If the cofinality of $(I,M^* \setminus I)$ is
$(\lambda^+,\lambda^+)$, then $I$ is definable in $M$ by the
property (e) of Claim A, so $B$ is definable in $M^*$.
So we may assume that the cofinality is not $(\lambda^+,\lambda^+)$.

First suppose that $ {\rm cf}(I) \leq \lambda$.
Then we can choose a set
$\{a_i : i< \lambda\}$ which  is cofinal in $I$.
Choose $j<\lambda^+$ with ${\rm cf}(j)=\lambda$ and 
$\{a_i : i< \lambda\} \subset M_j$.
If $M_j  \setminus I$ is bounded from below in
$M^* \setminus I$, say by $d \in M^* \setminus I$, then
$I$ is defined in $M^*$ by the formula
$\exists y [ x < y < d \wedge  y <_{\mathrm lev} e ]$,
where $e$ is an element from  $M_{j+1} \setminus M_j$.
So we may assume that there is a set $\{a_i':i<\lambda\}
\subset M_j \setminus I$  which is
coinitial in $M^* \setminus I$.
(We shall derive a contradiction from this. )
Let $b_{E_2} \in B$ with $b \notin M_j$.  Since the other
case can be treated similarly, we can
assume that $b \in I$. Then $b_{E_2}$ is
included in some interval $[0, a_i]$.
By the definition of $I$, there is $c_{E_2} \in B$ such
that $b_{E_2} <_{\mathrm lev} c_{E_2}$ and $a_i < c$.
But then $b$ and $c$ determine different Dedekind cuts
of $M_j$, hence $b$ and $c$ are not comparable
with respect to $<_{\mathrm tr}$. This contradicts our assumption
that $B$ is a branch.

Second suppose that the coinitiality of $M^* \setminus I$
is $\leq \lambda$ and that the cofinality of $I$ is $\lambda^+$.
As in the first case, we can choose $j<\lambda^+$ such that
$M_j \setminus I$ is coinitial in $M^* \setminus I$.
Choose $d \in I $ which bounds $I \cap M_j$ from above and
an element $e \in M_{j+1} \setminus M_j$.
Then $I$ is defined by the formula
$
\forall y[ d<y  \wedge y <_{\mathrm lev} e \rightarrow x < y].
$
Lastly the case where the cofinality of $(I, M^* \setminus I)$ is
$(\mu_1,\mu_2)$ with $\mu_1, \mu_2 \leq \lambda$ is
impossible by the definition of branch.
\qed

Let $T^*$ be the $L^*$-theory of $M^*$.
Under the hypothesis of Claim A
(i.e. $\diamondsuit_{S_\lambda^{\lambda^+}}$ etc),
we have proven the existence of $M^* \models T^*$
having a tree with the property
stated in Claim B.
However, by the absoluteness
(e.g. Thorem 6 in [1]),  the existence of such a model
can be proven without the hypothesis.
Moreover, as $T^*$ is countable, we can assume that relevant
properities of $M^*$ expressed by one
$L^*_{\omega_1\omega}(Q)$-sentence are
also possessed by such models.  ($Q$ is the quantifier which
expresses ``there are uncountably many''.)
Thus in ZFC we can show

\begin{claim}
There is a model $N^* \models T^*$ of cardinality $\aleph_1$
that satisfies:
\begin{enumerate}
\item The tree $(N^*/E_2, <_{\mathrm tr})$ has no undefinable branch;
\item The set $N^*/E_1$ of levels has the cardinality $\aleph_1$, but
for each $b/E_1 \in N^*/E_1$,
$\{c/E_1: c/E_1 <_{\mathrm lev}b/E_1\}$ is countable;
\item If $I$ is a definable subset of $N^*$ with
the Dedekind cut $(I,N^* \setminus I)$  of cofinality
$(\aleph_1,\aleph_1)$, then $I$ is definable in $N$;
\item The clause (d) of Claim A,
namely, for each level $d_{E_1}$ there is $a \in N^*$ such that
if  $b \in N^* \setminus \{S^n(0): n \in \omega\}$
then $\{F(a,c):  c < b\}$ includes  $\{c \in I: c \leq_{\mathrm lev} d\}$.

\end{enumerate}
\end{claim}

\begin{claim}
Let $N^*$ be a model of $T^*$ with the properties stated in
Claim C.
Then the reduct $N$ of $N^*$ to the language $L$ is
$\psi^*$-appropriate.
\end{claim}

\proof
Toward a contradiction, we assume that there is an undefinable
(in the sense of $N$) subset $I \subset N$
with $(N,I) \models \psi^*(P)$ and $I \neq \{S^n(0): n \in \omega\}$.
We show that the cofinality of $(I, N^* \setminus I)$ is
$(\aleph_1,\aleph_1)$.  Suppose that this is not the case.
First assume that the cofinality of $(I,<)$ is less than $\aleph_1$.
As $(N^*/E_1, <_{\mathrm lev})$ has the cofinality $\aleph_1$,
there is $d/E_1$ such that $\{c \in I: c \leq_{\mathrm lev} d\}$
is unbounded in $I$.  Since $I \neq \{S^n(0): n \in \omega\}$,
we can choose
$b \in I \setminus \{S^n(0):n \in \omega\}$.
By the fourth condition of Claim C, there is $a \in N^*$
such that
$\{F(a,c):  c < b\}$ includes  $\{c \in I: c \leq_{\mathrm lev} d\}$.
So $\{F(a,c):  c < b\} \cap I$ is unbounded in
$I$. This contradicts the last clause in the definition of
$\psi^*$.
Second assume that the coinitiality of $N^* \setminus I$ is
less than $\aleph_1$.
For a similar reason as in the first case,
we can find $d_{E_1}$ such that $\{c \in N^* \setminus I:
c \leq_{\mathrm lev} d\}$
is unbounded from below in $N^* \setminus I$.
Also we can choose $a \in N^*$ and $b \in I$ such that
$\{F(a,c):  c < b\}$ includes  $\{c \in I: c \leq_{\mathrm lev} d\}$.
If $I \cap \{F(a,c): c < b\}$ were bounded
(from above) say by $e \in I$,
then $I$ would be definable in $N$ by the $L$-formula
$$
\varphi(x,a,b,e) \stackrel{\rm def}{\equiv}
\forall z
[(e < z \wedge \exists y ( y<b \wedge z=F(a,y) ))\rightarrow x < z],
$$
contradicting our assumption that $I$ is not definable.
So $I \cap \{F(a,c): M_{i^*} \models c < b\}$ is not
bounded in $I$. Again this contradicts the last clause in
the definition of $\psi^*$.
So we have proven that the cofinality of $(I,N^* \setminus I)$
is $(\aleph_1,\aleph_1)$.

As in the proof of Claim B, we shall define a set
$\{(b_i)_{E_2} :i<\aleph_1\}$ and definable
intervals $J_i \subset N^*$ $(i<\aleph_1)$ such that
for each $i<\aleph_1$,
\begin{itemize}
\item $J_i$'s are decreasing;
\item $b_i \in J_i$, $J_i \cap I \neq \emptyset$, $J_i \cap (N^* \setminus I) \neq \emptyset$;
\item there is no element $d \in J_i$ with
$d <_{\mathrm lev} b_i$.
\end{itemize}
Suppose that we have chosen $d_j$'s and $J_j$'s for all $j < i$.
Since the cofinality of $I$ and the coinitiality of $N^* \setminus I$
are both $\aleph_1$, $\bigcap_{j<i} J_i$ intersects both  $I$ and
$N^* \setminus I$. Choose $b \in \bigcap_{j<i} J_i \cap I$ and
$c \in \bigcap_{j<i} J_i \cap (N^* \setminus I)$.
Then we put $J_i =\{e \in N^*: N^* \models b <e <d\}$.
Choose $b_i \in J_i$ of the minimum level.
(Such $b_i$  exists and $(b_i)_{E_2}$ is unique,
because every nonempty definable
subset of $N^*/E_1$ has  the minimum element with respect to
$<_{\mathrm lev}$. If there are
two such elements, they are distinguished by elements of
lower levels, contradicting the minimality.)
We claim that $\{(b_i)_{E_2}: i<\aleph_1\}$ determines a branch
$B=\{c_{E_2}: c_{E_2} \leq_{\mathrm tr} (b_i)_{E_2}
\text{ for some $i$}\}$.
For this it is sufficient to show that the $b_i$'s are linearly ordered by
$\leq_{\mathrm tr}$.
Let $i \leq i' <\aleph_1$. Then both $b_i$ and $b_{i'}$ are members of
the interval $J_i$.
Suppose that  $b_i$ and $b_{i'}$ are not comparable with respect to
$ \leq_{\mathrm tr}$.
They determine different Dedekind cuts of the elements
of lower levels.  So there is an element $c \in J_i$ with
$c <_{\mathrm lev} b_i$.  This contradicts our choice of
$b_i \in J_i$.
By our assumption (the fourth condition in Claim C),
the branch $B=\{(b_i)_{E_2}:i<\aleph_1\}$ is definable
in $N^*$.
It is easy to see that $I$ and $B$ are interdefinable in $N^*$.
So $I$ is also definable in $N^*$, hence $I$ is definable in $N$ by the
third condition in Claim C.
This contradicts our assumption that $I$ is undefinable in $N$.
\qed


\begin{center}
{\bf References.}
\end{center}

[1] S. Shelah,
 Models with second-order properties II.
 Trees with no undefined branches, Annals of Mathematical
 Logic 14 (1978), pp. 73-87. [Sh:73].

[2] S. Shelah,
Models with second order properties IV.
A general method and eliminating diamonds, Annals of Pure and
Applied Logic 25 (1983), pp. 183-212.
[Sh:107].

[3] S. Shelah,
Non structure theory. In preparation. [Sh:e].

\end{document}